\documentclass[11pt]{article}
\usepackage{amsfonts,latexsym}
\usepackage{fullpage}
\usepackage{amscd}
\usepackage[all]{xy}
\usepackage{float,amsmath,amssymb,mathrsfs,bm,multirow,graphics}
\usepackage{graphicx}
\usepackage{color}

\newcommand{\proofbegin}{\noindent{\it Proof.\,\,}}
\newcommand{\proofend}{\hfill$\Box$\bigskip}
\newcommand{\be}{\begin{equation}}
\newcommand{\ee}{\end{equation}}

\newcommand{\R}{\mathbb{R}}
\newcommand{\Z}{\mathbb{Z}}
\newcommand{\T}{\mathbb{T}}
\newcommand{\diag}{\mbox{diag}}
\newcommand{\Je}{J_{\epsilon}}
\newcommand{\Ue}{U_{\epsilon}}

\newtheorem{theorem}{Theorem}

\newtheorem{lemma}[theorem]{Lemma}

\title{Nonuniform dependence on initial data for compressible gas dynamics:
The periodic Cauchy problem}

\author{Barbara L. Keyfitz and Feride T\i\u{g}lay}

\begin{document}
\maketitle
\bibliographystyle{abbrv}
 \begin{abstract}
 We start with the classic result that the Cauchy problem for ideal compressible gas dynamics is locally well posed in time in the sense of Hadamard; there is a unique solution that depends continuously on initial data in Sobolev space $H^s$ for $s>d/2+1$ where $d$  is the space dimension. We prove that the data to solution map for periodic data in two dimensions although continuous is not uniformly continuous on any bounded subset of Sobolev class functions.
 \end{abstract}

The compressible gas dynamics equations of ideal hydrodynamics
are given by the system
\begin{equation} \label{eq:ncompgas1}
\begin{split}
\rho_t +\nabla\cdot(\rho {\bf u})&=0 \\
(\rho {\bf u})_t+\nabla\left(\rho {\bf u} \times{\bf u} \right)+\nabla p&=0\\
(E\rho)_t+ \nabla\cdot\left(E\rho {\bf u}  + p{\bf u} \right)&=0
\end{split}
\end{equation}
with ${\displaystyle E=
e+\frac{1}{2}|{\bf u}|^2}$ the total energy
and ${\displaystyle  e=\frac{p}{(\gamma -1)\rho}}$ the
internal energy, expressed in terms of density $\rho$,
pressure $p$ and velocity $\bf u$.

Classical solutions and well-posedness in Sobolev spaces
(existence and uniqueness of solutions as well as continuous 
dependence of solutions on initial data)
of the initial value problem for  
\eqref{eq:ncompgas1} have been studied extensively, 
see for instance \cite{Kato}, \cite{Majda1}, \cite{Serre}  and \cite{Taylor}. 
Sobolev space results are all local in time. 
In one space dimension shock waves form in finite time for almost all data in $H^s$,
and for later times only weak solutions exist.
(The definition of weak solutions, and well-posedness theory in
$BV_{\textrm{loc}}\cap L^1_{\textrm{loc}}$, which are not the
subject of this paper, can be found in \cite{Bre:book} and \cite{Dafbook}.) 
In higher dimensions there is as yet no existence theory for weak solutions,
and classical (Sobolev space) solutions have a finite-time life span for almost
all data  \cite{Majda1,Taylor}.

Our goal is to study continuity properties of the solution map for classical
solutions; in this paper
we prove that for periodic data the initial-data to solution map is not uniformly continuous 
in Sobolev spaces.
In a  companion paper, \cite{HoKeTi}, we extend this result to $H^s$ data
in the plane.
Throughout, we assume $s$ to be large enough for classical results to hold.

We consider solutions $U=U({\bf x},t)$ that take values in 
a compact subset of the state space $G=\{U\equiv(\rho, {\bf u}, p)\mid \rho, p>0\}$,
defined as the region where the
physical quantities  $\rho$ and  $e$ are positive, 
and the system is symmetrizable hyperbolic.  

In two dimensions, since we are considering classical solutions, 
we can ignore conservation form and 
write system \eqref{eq:ncompgas1} as
\begin{equation}\label{eq:classical1}
\begin{split}
 \rho_t +u\rho_x+v\rho_y+\rho(u_x+v_y)&=0 \\
 u_t +uu_x+ vu_y+ h_x+ \frac{h}{\rho}\rho_x&=0 \\
v_t +uv_x+ vv_y+ h_y+  \frac{h}{\rho}\rho_y&=0 \\
h_t+uh_x+ vh_y+(\gamma -1) h(u_x+v_y)&=0\,.
\end{split}
\end{equation}
The parameter $\gamma$ denotes the ratio of specific heats (typically $1<\gamma
<3$) 
and $h=p/\rho = (\gamma-1)e$ is a multiple of the internal energy.  

We study this system in Sobolev spaces on the two dimensional torus: 
$H^s(\T^2)$ where $\T=\R/2\pi\Z$. The Sobolev norm is given by
\[ \|u\|_{s}^2=\langle \Lambda^s u, \Lambda^s u \rangle\,,
\]
where $\Lambda^s=(1-\Delta)^{s/2}$ and $\langle \ , \ \rangle$ denotes the 
$L^2$ inner product. 
Defining $U=(\rho,u,v,h)$ and $U(t)=U(\cdot,t)$,
our main result is
\begin{theorem}[Nonuniform Dependence on Initial Data] \label{th:nonuni}
For $s>2$, the data to solution map $U(0)\to U(t)$ for the system \eqref{eq:classical1}  is not 
uniformly continuous from any bounded subset of $(H^s(\T^2))^4$ into 
$C([0,T]; (H^s(\T^2))^4)$.
\end{theorem}
We note the significance of $s>2$.
The well-posedness theory for symmetrizable hyperbolic systems, 
which forms the basis for our analysis,
is credited to G\aa rding \cite{ga}, Leray \cite{leray}, Kato \cite{Kato} and Lax \cite{Lax4}.
Solutions for quasilinear systems in $d$ space dimensions exist in spaces $H^s$
for $s> d/2 +1$.
Modern expositions of the theory can be found in Majda \cite{Maj2}, Serre \cite{Serre}
or Taylor \cite{Taylor3}.

We give the proof of Theorem \ref{th:nonuni} in Section \ref{thmpf}.
Our proof uses a framework introduced to prove an analogous result for 
the incompressible Euler equations of ideal hydrodynamics in \cite{Himonas-Misiolek}. 
This framework has been used for other nonlinear PDE including the Benjamin-Ono 
equation in \cite{KochTzvet} and the Camassa-Holm equation on the real line and on the one 
dimensional torus in \cite{HimKen} and  \cite{HimKenMis} respectively. 
Implementation of this framework for the periodic Cauchy 
problem for the incompressible two-dimensional Euler equations
is carried out in \cite{Himonas-Misiolek} with minimal technicalities.
In that case, 
two sequences of exact solutions $\{U_{-1,n}(t)\}$ and $\{U_{1,n}(t)\}$ in $H^s$ are 
constructed such that
as  $n\rightarrow \infty$, 
\begin{equation}
\label{eq:nonuniexact}
 \| U_{-1,n}(0)-U_{1,n}(0)\|_s \rightarrow 0 \mbox{ and } \| U_{-1,n}(t)-U_{1,n}(t)\|_s 
 \geq \sin t \mbox{ for } t>0\,.
\end{equation}
Exact solutions with this property exist for the compressible system as
well, as we show in Section \ref{cstden}, but they have the unsatisfactory
feature of being almost trivial:
They have constant density and pressure (they are thus also solutions
of the incompressible equations).
Our proof of Theorem \ref{th:nonuni} exhibits the phenomenon of nonuniform 
dependence in a situation
where density and pressure also vary, by adapting the Himonas-Misio\l ek
construction in  \cite{Himonas-Misiolek}.
As exact solutions of \eqref{eq:classical1} with non-constant density are not available,  
we use instead two sequences, similar to those constructed in \cite{Himonas-Misiolek},
which we prove are approximate solutions.
Section \ref{wplf} sets up the background for the construction, and
in Section \ref{nondep} we prove the critical estimate  that shows the
approximate solutions are close enough to 
exact solutions to give the estimates \eqref{eq:nonuniexact} for actual solutions.
The final section, Section \ref{conc}, includes some comments on the examples
and on the significance of the result.
\section{Well-posedness and Lifespan} \label{wplf}
In this section, we present a suggestive example, and
review some of the classical results, mentioned in
the introduction, for system \eqref{eq:ncompgas1} or \eqref{eq:classical1}.

\subsection{A Constant-Density Example}
\label{cstden}
The following example presents a pair
of sequences, somewhat simpler than the exact solutions of 
\cite{Himonas-Misiolek},
that solve both the incompressible and the compressible gas dynamics
system, and are easily seen to have the property \eqref{eq:nonuniexact}.
The functions
\be \label{Vseq} 
V_{\omega,n}(x,y,t)= (\rho,u,v,h)=
\left(\rho_0,\frac{1}{n^s}\cos(ny-\omega t), \frac{\omega}{n},h_0\right)\ee
for $\omega = \pm 1$
are exact solutions of \eqref{eq:classical1}.
Each solution is divergence-free;
in Section \ref{conc} we note that these sequences also satisfy the
incompressible system, \eqref{incomp}.
(Solutions of this form may be known but it seems not to have been observed that they exhibit this property.)
We carry out verification of \eqref{eq:nonuniexact}, which is straightforward.
For each $n$, 
$$V_{1,n}(x,y,0)-V_{-1,n}(x,y,0)=\left(0,0,\frac{2}{n},0\right)\,,$$
and clearly this tends to $0$ in $H^\sigma$ for any $\sigma \geq 0$.
On the other hand,
$$V_{1,n}(x,y,t)-V_{-1,n}(x,y,t)=\left(0,\frac{2}{n^s}\sin ny\sin t,\frac{2}{n},0\right)\,.$$
A straightforward calculation (see \cite[Lemma 3.2]{Himonas-Misiolek}) gives the values
\begin{equation} \label{Hsnorm}
\| \cos nx\|_{\sigma}= \| \sin nx\|_\sigma=
\pi\sqrt{2}(1+n^2)^{\sigma/2}\end{equation}
for the one-dimensional $H^\sigma(\T)$ norms for any $\sigma$, and so
\be \label{goodest}
\|V_{1,n}(\cdot,t)-V_{-1,n}(\cdot,t)\|_s = \frac{2\sqrt{2}\pi}{n^s}(1+n^2)^{s/2}
|\sin t| + \frac{2}{n} \gtrsim |\sin t|\,;\ee
that is, the difference in $H^s$ between two solutions does not go to zero
for $t\neq 0$.
(The notation $\lesssim$, $\gtrsim$ and $\simeq$ indicates that the relations hold
up to constants independent of $n$.)

The approximate solutions we construct for our proof of Theorem
\ref{th:nonuni} exhibit non-uniform dependence
on data via the same mechanism.
Their structure is similar to, but not quite the same as, the solutions \eqref{Vseq}.
We emphasize that the actual solutions to \eqref{eq:classical1} with the same
initial data as the approximate solutions \eqref{eq:happ} below do not have constant
density.
In particular, they all develop shocks, but after a time that is bounded away from
zero, uniformly in $n$.

This example, simple as it is,
forms the basis for the demonstration of
non-uniform dependence in $H^s(\mathbb R^2)$,
both for Himonas and Misio\l ek in \cite{Himonas-Misiolek}
and for our adaptation for the compressible equations \eqref{eq:classical1} in 
a companion paper, \cite{HoKeTi}.
When transforming periodic data to $H^s$-integrable data by introducing  cut-off
functions, one introduces perturbations to the density and pressure,
so the full-plane variant of this example is not a constant-density
solution.

\subsection{Symmetrized System}
The equations for compressible ideal gas dynamics \eqref{eq:ncompgas1}  
form a classical model from mathematical physics, one that indeed
motivated the theory of symmetric and symmetrizable hyperbolic systems. 
We express  system \eqref{eq:classical1} in the form
\begin{equation} \label{classical}
U_t+A(U)U_x+B(U)U_y=0
\end{equation}
with 
\[ U=  \begin{pmatrix}
\rho \\ 
u\\
v\\
h
\end{pmatrix}\,,  \qquad 
 A(U)= \begin{pmatrix}
u & \rho & 0 & 0 \\ 
h/\rho & u & 0 & 1 \\
0 & 0 & u & 0 \\
0 &( \gamma-1) h & 0 & u
\end{pmatrix}\,, \qquad  
B(U)=   \begin{pmatrix}
v & 0 & \rho & 0 \\ 
0 & v & 0 & 0 \\
h/\rho & 0 & v & 1 \\
0 & 0 & (\gamma-1) h & v
\end{pmatrix} \,, \]
and note that it is symmetrizable.
If we let 
\[ A_0(U) =  \begin{pmatrix}
h/\rho & 0 & 0 & 0 \\ 
0 & \rho & 0 & 0 \\
0 & 0 & \rho & 0 \\
0 &0 & 0 & \frac{\rho}{(\gamma -1)h}
\end{pmatrix} \,,
\]
then $A_0(U)$ is a positive definite symmetric matrix for $U\in G$
and we have the equivalent symmetric hyperbolic system
\[ A_0 U_t+A_1(U)U_x+B_1(U)U_y= 0
\]
with 
\[ 
 A_1(U)=\begin{pmatrix}
\frac{uh}{\rho} & h & 0 & 0 \\ 
h & \rho u& 0 & \rho \\
0 & 0 & \rho u& 0 \\
0 & \rho & 0 & \frac{\rho u}{( \gamma-1) h} 
\end{pmatrix}\,, \qquad  
B_1(U)=  \begin{pmatrix}
\frac{v h}{\rho} & 0 & h & 0 \\ 
0 & \rho v & 0 & 0 \\
h & 0 & \rho v & \rho \\
0 & 0 & \rho & \frac{\rho v}{(\gamma-1) h}
\end{pmatrix}\,. \]

\subsection{Lifespan and Solution Size Estimates}
A standard approach in proving existence and uniqueness of solutions for 
Cauchy problems is to obtain a solution as a limit to a mollified system.
This is the approach taken by Taylor, \cite{Taylor3}. 
Let $U_{\epsilon}$ be a solution of 
\begin{equation}
A_0(\Je\Ue) \partial_{t}\Ue + A_1(\Je\Ue)\partial_x (\Je\Ue)+
 B_1(\Je\Ue)\partial_y (\Je\Ue)=0
\end{equation}
where $\Je, \ 0< \epsilon \leq 1$ is a Friedrichs mollifier,
defined  by a Fourier series representation
\[ (\Je v)^{\wedge}(l)= \varphi(\epsilon l) \widehat{v}(l), \quad l\in \Z^2
\]
with $\varphi\in C_0^{\infty}(\R^2)$  real-valued and $\varphi(0)=1$.
Then the existence and uniqueness of solutions follow from a general 
argument for symmetrizable hyperbolic systems.  
The proof uses an energy estimate (see 
Chapter 16 in \cite{Taylor3} for instance,
or  estimate (2.50) in the statement of Theorem 2.2 in \cite{Maj2}) 
that leads to a solution size estimate
\begin{equation} \label{eq:solsize}
\|U(t)\|_s \leq C \|U_0\|_s \quad {\mbox for } \quad t\in [0,T]
\end{equation} 
where $T$ depends on $\|U_0\|_s$. 

Continuous dependence on initial data is shown by Kato \cite[Theorem III(b)]{Kato}, 
who proves
that there exist $d>0$ and $T'\in(0, T]$ depending only on 
the $H^s$ norm of the initial data 
$U_0$ such that if $\| U_0^n-U_0\|_{s}\leq d$ and $\lim_{n\rightarrow \infty}\| 
U_0^n-U_0\|_s=0$ then the solutions $U^n$ exist on a common interval 
$[0,T']$ and $\| U^n(t)-U(t)\|_{s}\rightarrow 0$ uniformly in $t$. 
The solution map $U_0 \rightarrow U$ is not H\"{o}lder 
continuous in the $H^s$ norm. 
The comparison between our result and Kato's is discussed in the
final section of this paper.

\section{Nonuniform Dependence} \label{nondep}

In this section we construct a set of approximate solutions, show that
they are good approximations to a true solution, and 
prove a critical estimate, Theorem \ref{prop:error}.
\subsection{Approximate Solutions}

Our strategy is to use two sequences 
$U^{\omega,n}= \left( \rho^{\omega,n}, u^{\omega,n}, v^{\omega,n}, 
h^{\omega,n} \right)^\intercal$, with
$\omega=\pm 1$, of approximate solutions:
\begin{equation}  \label{eq:happ}
\begin{split}
\rho^{\omega,n} &= \rho_0 \\
u^{\omega,n}&=\frac{\omega}{n}+\frac{1}{n^s}\cos(ny-\omega t) \\
v^{\omega,n}&=\frac{\omega}{n}+\frac{1}{n^s}\cos(nx-\omega t) \\
 h^{\omega,n}&=h_0+\frac{1}{n^{2s}}\sin(nx-\omega t) \sin(ny-\omega t)
 \end{split}
\end{equation}
that are arbitrarily close at time zero but are separated at later times. 
The approximate solutions are in $(H^{s})^4$ and their $H^s$ norms
 are uniformly bounded in $n$.

Let $U\equiv U_{\omega,n}$ represent the actual solution to 
\eqref{classical} with the same initial values as $U^{\omega,n}$:
\be \label{eq:initial}
 U_{\omega,n}(0)=U^{\omega,n}(0)= \left(
\rho_0, \;
\frac{\omega}{n}+\frac{1}{n^s}\cos ny, \;
\frac{\omega}{n}+\frac{1}{n^s}\cos nx, \;
h_0+\frac{1}{n^{2s}}\sin nx \sin ny
\right)^\intercal\,.
\ee
To estimate dependence of the solution size on $n$
we introduce the notation $\widetilde{U}\equiv(\tilde \rho,u,v,\tilde h)^\intercal
=(\rho-\rho_0, u, v, h-h_0)^\intercal$, subtracting the stationary solution $(\rho_0,0,0,h_0)$
from both the approximate and the actual solutions.

From \eqref{Hsnorm} we have, for any $\sigma \geq 0$,
\begin{equation} \label{eq:appest} 
\|\widetilde{U}^{\omega,n}\|_{\sigma} \leq C n^{\sigma-s}\,.
\end{equation}
The solution size estimate \eqref{eq:solsize} also applies to functions $\widetilde U$
derived from the exact solutions to \eqref{classical}, since $\widetilde U$ satisfies \eqref{classical}
with modified but still symmetrizable coefficients, so the same estimates from \cite{Taylor3}
give us \eqref{eq:solsize} and thence \eqref{eq:appest} for 
$\widetilde U_{\omega,n}=U_{\omega,n}-(\rho_0,0,0,h_0)$.

Another calculation shows that the approximate solutions satisfy the equation
\[ U^{\omega,n}_t+A(U^{\omega,n})U^{\omega,n}_x+B(U^{\omega,n})
U^{\omega,n}_y= (0, 0, 0, R_4)^\intercal\,,
\]
where the residue is given by
\begin{align*} R_4&= \frac{1}{n^{3s-1}}\cos(nx-\omega t) \cos(ny-\omega t)\big( 
 \sin(nx-\omega t)+\sin(ny-\omega t) \big).\\
 &= \frac{1}{2n^{3s-1}}\big(\sin 2(nx-\omega t)\cos (ny-\omega t)
 + \cos(nx-\omega t)\sin 2(ny-\omega t)\big)\,.
\end{align*}

\begin{lemma}[Residue Estimate]  For $n\gg 1$, $1<\sigma\leq s-1$ and
 $s>2$ the residue satisfies
\[ \| R_4 \|_{\sigma} \leq C n^{2\sigma -3s+1}.
\]
\label{res_est}
\end{lemma}
\proofbegin
The estimate follows from  the one-dimensional norms, \eqref{Hsnorm}.
\proofend

\subsection{Error Estimates}

We fix $\omega$ and $n$ and let $U$ and $\widetilde U$
denote $U_{\omega,n}$ and $\widetilde U_{\omega,n}$.
Our goal in this section is to calculate the error $\mathcal{E}=U - U^{\omega,n}
\equiv(E, F, G, H)^\intercal$, the difference between actual and approximate solutions,
and show that it goes to zero in the $H^s$ norm as $n\rightarrow \infty$. 
The error $\mathcal{E}$ satisfies the system of equations
\begin{equation} 
\label{eq:errsys}
 \mathcal{E}_t+A(U^{\omega,n})\mathcal{E}_x+B(U^{\omega,n})\mathcal{E}_y+
 C(U)\mathcal{E}+ (0, 0, 0, R_4)^\intercal=0\,,
\end{equation}
where 
\[ C(U)=  \begin{pmatrix}
u_x+v_y & \rho_x & \rho_y & 0 \\ \\
-\frac{h^{\omega,n}\rho_x}{\rho \rho_0} & u_x & u_y & \frac{\rho_x}{\rho} \\ \\
-\frac{h^{\omega,n}\rho_y}{\rho \rho_0} & v_x & v_y & \frac{\rho_y}{\rho} \\ \\
0 & h_x & h_y & (\gamma-1)(u_x+v_y)
\end{pmatrix}\,. \]
To obtain the desired estimates, we work
in a second Sobolev space, $H^\sigma$, with  $1<\sigma < s-1$.
One of the tools we use is the following commutator estimate, 
which is a special case of Proposition 4.2 from \cite{Taylor}:
\begin{lemma}[Commutator Lemma] \label{commut}
For $k>2$ and $1<\sigma \leq k$,
\begin{equation} \label{eq:est_com}
\| [\Lambda^{\sigma},f]u\|_{L^2}\leq C \| f\|_{k}\|u\|_{\sigma -1}\,,
\end{equation}
where $ [\Lambda^{\sigma},f]u = \Lambda^{\sigma} (fu) - f \Lambda^{\sigma} u$.
\end{lemma}
We also need the following lemma.
\begin{lemma}[Reciprocal Lemma] \label{lem:rec}
For $s>1$ and $\sigma\leq s$ let $f\in H^\sigma(\T^2)$ and suppose the density 
$\rho\in H^s(\T^2)$ is in a compact 
subset of the state space $G$. 
Then $f/\rho\in H^\sigma(\T^2)$ and 
\begin{equation}
\left\| \frac{f}{\rho} \right\|_\sigma \leq C \left(1+ \| \tilde{\rho} \|_s^\sigma\right) \|f \|_\sigma.
\end{equation}
\end{lemma}
The proof of this lemma is given in \cite{Kato} 
(Lemma 2.13 and the argument following) for integer values of $s$ and $\sigma$. 
For the non-integer case, a proof is given in \cite{HoKeTi}.

The approximate solutions exhibit non-uniform dependence
via an argument, given in Section \ref{thmpf},
similar to that presented in Section \ref{cstden}. 
Thus, the
heart of the nonuniform dependence theorem, Theorem \ref{th:nonuni},
is the demonstration that the approximations are indeed $H^s$-close
to an actual solution.
The crucial technical estimate is the following theorem.
It is established in a Sobolev space with index strictly smaller
than the space of interest.
We will see that this suffices.

\begin{theorem} \label{prop:error}
The system \eqref{eq:errsys} is symmetrizable
and for $s>2$, $1<\sigma < s-1$ and $n\gg 1$ the error 
$\mathcal{E}=U - U^{\omega,n}$ satisfies the estimate
\begin{align}
 \| \mathcal{E}(t) \|_{\sigma} & \leq  n^\beta \left( e^{ct}-1 \right)\,, \quad
 \textrm{where}\quad \beta = \max\{2\sigma-3s+2, \sigma-2s\}\,,
 \label{eq:errm}
 \end{align}
and  $c$ depends on 
 $\rho_0$, $h_0$ and $\gamma$ and decreases with $n$.
\end{theorem}

\proofbegin
Upon multiplying  the system \eqref{eq:errsys} by the 
symmetric matrix $ A_0(U^{\omega, n})$,  
the symmetrized system for the error is
\begin{equation} 
\label{eq:errsyssym}
A_0(U^{\omega, n})\mathcal{E}_t+A_1(U^{\omega, n})\mathcal{E}_x+
B_1(U^{\omega, n})\mathcal{E}_y+C_1(U^{\omega, n},U)\mathcal{E}+ 
A_0(U^{\omega, n})(0, 0, 0, R_4)^\intercal=0\,,
\end{equation}
where $C_1(U^{\omega, n},U)=A_0(U^{\omega, n})C(U)$.

We apply  $\Lambda^{\sigma}$ to \eqref{eq:errsyssym} and take the
$L^2$ inner product with $\Lambda^\sigma\mathcal E$ to
 obtain 
\begin{align}
\langle \Lambda^{\sigma} \mathcal{E}, \Lambda^{\sigma} \left( A_0(U^{\omega,n})
\mathcal{E}_t \right) \rangle = & - \langle \Lambda^{\sigma} \mathcal{E}, 
\Lambda^{\sigma} \left( C_1(U^{\omega,n},U)\mathcal{E} \right) \rangle \label{eq:gr1}\\
& - \langle \Lambda^{\sigma} \mathcal{E} , \Lambda^{\sigma} 
\left( \diag(A_1(U^{\omega,n}))\mathcal{E}_x + \diag(B_1(U^{\omega,n}))\mathcal{E}_y 
\right)\rangle \label{eq:gr2}\\
& -\langle \Lambda^{\sigma} \mathcal{E}, \Lambda^{\sigma} \left( A_R(U^{\omega,n})
\mathcal{E}_x + B_R(U^{\omega,n})\mathcal{E}_y \right)\rangle \label{eq:gr3}\\
& -  \left\langle \Lambda^{\sigma} H , \Lambda^{\sigma} \left( \frac{\rho_0}
{(\gamma-1)h^{\omega,n}} R_4 \right)\right\rangle \label{eq:gr4}\,,
\end{align}
where $\diag(A)$ denotes the diagonal part of a matrix $A$ and $A_R=A - \diag(A)$. 

The first step is to establish the estimate
\begin{equation}
\big| \langle \Lambda^{\sigma} \mathcal{E}, \Lambda^{\sigma}\big( A_0(U^{\omega,n})
\mathcal{E}_t\big) \rangle\big| \leq C \big[ 
n^{\max\{-1,\sigma-s+1\}}  \| \mathcal{E}\|_{\sigma}^2 +
n^{2\sigma-3s+1} \|\mathcal{E}\|_{\sigma}\big]\,,
\label{eq:A0Et}
\end{equation}
where $C$ depends only on $\rho_0$, $\gamma$ and $h_0$.

With a change of sign, the first expression, \eqref{eq:gr1}, is 
\begin{multline} \label{gr1expand}
\left\langle \Lambda^{\sigma} E ,  \Lambda^{\sigma} \left(\frac{h^{\omega,n}}{\rho_0}
(u_x+v_y)E+\frac{h^{\omega,n}}{\rho_0}\rho_xF+
\frac{h^{\omega,n}}{\rho_0}\rho_y G\right)\right\rangle \\+
\left\langle \Lambda^{\sigma} F ,  \Lambda^{\sigma} \left(-\frac{h^{\omega,n}\rho_x}
{\rho}E+u_x\rho_0F+u_y\rho_0G+\frac{\rho_0\rho_x}{\rho}H\right)\right\rangle \\+
\left\langle \Lambda^{\sigma} G ,  \Lambda^{\sigma} \left(-\frac{h^{\omega,n}\rho_y}
{\rho}E+ \rho_0v_xF+ \rho_0v_yG+\frac{\rho_0\rho_y}{\rho}H \right)\right\rangle \\+
\left\langle \Lambda^{\sigma} H ,  \Lambda^{\sigma} \left( \frac{\rho_0h_x}
{(\gamma -1)h^{\omega,n}}F+ \frac{\rho_0h_y}{(\gamma -1)h^{\omega,n}}G
+\frac{\rho_0(u_x+v_y)}{(\gamma -1)h^{\omega,n}}H \right)\right\rangle\,.
\end{multline}
We use Cauchy-Schwarz on the first term in \eqref{gr1expand}:
\begin{multline*} T_1\equiv \left|
\left\langle \Lambda^{\sigma} E ,  \Lambda^{\sigma} \left(\frac{h^{\omega,n}}{\rho_0}
(u_x+v_y)E+\frac{h^{\omega,n}}{\rho_0}\rho_xF+
\frac{h^{\omega,n}}{\rho_0}\rho_y G\right)\right\rangle\right| \\ \leq C
 \|E\|_{\sigma} \left\| (h_0+\tilde{h}^{\omega,n}) \big((u_x+v_y)E+ 
 \rho_xF+\rho_yG\big)\right\|_{\sigma}\,,
 \end{multline*}
where $C$ depends on $\rho_0$. 
From the algebra property of Sobolev spaces \cite[page 106]{Adams}, valid for
 $\sigma > 1$, we obtain
\[ T_1 \leq C \|\mathcal{E}\|_{\sigma}^2 
( \| \widetilde{U}^{\omega,n} \|_{\sigma}+1) \| \widetilde{U}\|_{\sigma+1}. \]
By the solution size estimate \eqref{eq:solsize}, 
and the bound \eqref{eq:appest} applied to the initial data, 
$\|\widetilde{U}\|_{\sigma+1}$ is bounded, 
up to a constant independent of $n$, by 
$n^{\sigma+1-s}$. 
Using the same bound \eqref{eq:appest} for $\|\widetilde{U}^{\omega,n}\|_\sigma$ and 
noting that $2\sigma-2s+1<\sigma-s+1$, we obtain
\begin{equation}\label{eq:gr1exp1}
T_1 \leq C n^{\sigma-s+1} 
\|\mathcal{E}\|_{\sigma}^2\,.
\end{equation}

To estimate the second term in \eqref{gr1expand} we use Cauchy-Schwarz 
and the algebra property of Sobolev spaces as above to obtain
\begin{multline*}T_2\equiv\left| \left\langle \Lambda^{\sigma} F ,  
\Lambda^{\sigma} \left(-\frac{h^{\omega,n}\rho_x}
{\rho}E+u_x\rho_0F+u_y\rho_0G+\frac{\rho_0\rho_x}{\rho}H\right)\right\rangle\right| \\
\leq C \|F\|_{\sigma} \left( \|h^{\omega,n}\|_{\sigma} \left\|\frac{\rho_x}{\rho}\right\|_{\sigma}
 \|E\|_{\sigma} + \|u_x \|_{\sigma}  \| F\|_{\sigma} + \|u_y \|_{\sigma}  \| G\|_{\sigma}  +
 \left\|\frac{ \rho_x}{\rho}\right\|_{\sigma} \| H\|_{\sigma}\right) \\
\leq C \|\mathcal{E}\|_{\sigma}^2   \left( (\|\widetilde{U}^{\omega,n}\|_{\sigma}+1) 
\left\|\frac{\rho_x}{\rho}\right\|_{\sigma} + \|\widetilde{U} \|_{\sigma+1}\right).
\end{multline*}
Using the Reciprocal Lemma, Lemma \ref{lem:rec}, with the solution size 
estimate \eqref{eq:solsize},  applied to the derived
solution $\widetilde U$,
and the bound \eqref{eq:appest} applied to the initial data leads to
\[ T_2
\leq C  \|\mathcal{E}\|_{\sigma}^2  \left( 1+\|\widetilde{U}^{\omega,n}\|_{\sigma}\right) 
\left(1+  \|\widetilde{U} \|_{s}^{\sigma}\right) \|\widetilde{U}\|_{\sigma+1}\,. 
\]
Since $\sigma-s<0$ and $n\gg 1$, the largest power of $n$ in
this expression is $\sigma-n+1$; therefore
\begin{equation}\label{eq:gr1exp2}
T_2 \leq 
C n^{\sigma-s+1} \|\mathcal{E}\|_{\sigma}^2\,.
\end{equation}

The third term in \eqref{gr1expand} is estimated like the second term above and 
yields the same bound.

For the last term in  \eqref{gr1expand} we have the following estimate by Cauchy-Schwarz 
and the algebra property of Sobolev spaces:
\begin{multline*} T_3\equiv
\left|\left\langle \Lambda^{\sigma} H ,  \Lambda^{\sigma} \left( \frac{\rho_0h_x}
{(\gamma -1)h^{\omega,n}}F+ \frac{\rho_0h_y}{(\gamma -1)h^{\omega,n}}G
+\frac{\rho_0(u_x+v_y)}{(\gamma -1)h^{\omega,n}}H \right)\right\rangle\right| \\
\leq C \| H\|_{\sigma}  \left( \left\| \frac{h_x}{h^{\omega,n}}\right\|_{\sigma}  
\| F\|_{\sigma} + \left\| \frac{h_y}{h^{\omega,n}}\right\|_{\sigma} 
\| G\|_{\sigma} + \left\| \frac{u_x+v_y}{h^{\omega,n}}\right\|_{\sigma} \| H\|_{\sigma}  \right) 
\end{multline*}
where $C$ depends on $\rho_0$ and $\gamma$. 
Using the Reciprocal Lemma \ref{lem:rec} with the bound \eqref{eq:appest} 
applied to the initial data and the solution size estimate \eqref{eq:solsize} leads to
\begin{equation} \label{eq:gr1exp4}
T_3 \leq C n^{\sigma-s+1} \|\mathcal{E} \|_{\sigma}^2.
\end{equation}
Combining the estimates \eqref{eq:gr1exp1} - \eqref{eq:gr1exp4} 
we obtain a bound for \eqref{eq:gr1}:
\begin{equation}\label{gr1est}
\left|\langle \Lambda^{\sigma} \mathcal{E}, 
\Lambda^{\sigma} \left( C_1(U^{\omega,n},U)\mathcal{E} \right) \rangle \right|
\leq C n^{\sigma-s+1} \|\mathcal{E} \|_{\sigma}^2.
\end{equation}

The  expression \eqref{eq:gr2}, with a change of sign, is
\begin{multline*}
\left\langle \Lambda^{\sigma} E ,\Lambda^{\sigma} \left(
\frac{h^{\omega,n}u^{\omega,n}}{\rho_0} E_x  + \frac{h^{\omega,n}
 v^{\omega,n}}{\rho_0} E_y \right)\right\rangle +
\langle \Lambda^{\sigma} F , \Lambda^{\sigma} \left(\rho_0u^{\omega,n}F_x 
+\rho_0v^{\omega,n}F_y \right)\rangle \\
+\langle \Lambda^{\sigma} G , 
\Lambda^{\sigma} \left(\rho_0 u^{\omega,n}G_x + \rho_0v^{\omega,n}G_y
\right)\rangle +
\left\langle \Lambda^{\sigma} H ,  \Lambda^{\sigma}\left(\frac{\rho_0u^{\omega,n}}
{(\gamma-1)h^{\omega,n}}H_x + \frac{\rho_0v^{\omega,n}}{(\gamma-1)h^{\omega,n}}
H_y \right)\right\rangle\,.
\end{multline*}
All terms are estimated in the same way; we demonstrate the details of the first
by writing $\langle \Lambda^{\sigma} E ,  \Lambda^{\sigma} \left(
\frac{h^{\omega,n}u^{\omega,n}}{\rho_0} E_x \right)\rangle$ using commutators:
\begin{eqnarray}
&\langle \Lambda^{\sigma} E ,  \Lambda^{\sigma} \left(\frac{h^{\omega,n}
u^{\omega,n}}{\rho_0} E_x \right)\rangle = &\langle \Lambda^{\sigma} E ,  
\left[\Lambda^{\sigma},  \frac{h^{\omega,n}u^{\omega,n}}{\rho_0} \right] 
E_x\rangle \label{eq:sec1} \\
&& + \langle \Lambda^{\sigma} E , \frac{h^{\omega,n}u^{\omega,n}}
{\rho_0}\Lambda^{\sigma} E_x \rangle \label{eq:sec2} 
\end{eqnarray}
Using the commutator estimate \eqref{eq:est_com} with $k=\sigma+1$ in 
\eqref{eq:sec1}  and taking account of \eqref{eq:appest} we have
\begin{multline*}
\left|\langle \Lambda^{\sigma} E ,  \left[\Lambda^{\sigma},  \frac{h^{\omega,n}u^{\omega,n}}
{\rho_0} \right] E_x\rangle\right| \leq C \| (h_0+{\tilde h}^{\omega,n}) 
u^{\omega,n}\|_{\sigma+1}\|E\|_{\sigma}^2 
 \leq C  (1+\|{\widetilde U}^{\omega,n}\|_{\sigma+1}) \|\widetilde{U}^{\omega,n}
 \|_{\sigma+1}\|E\|_{\sigma}^2 \\
  \leq C n^{\max\{2(\sigma-s+1), \sigma-s+1\} } \|E\|_{\sigma}^2 \leq C n^{\sigma-s+1} 
  \|E\|_{\sigma}^2.
\end{multline*}
We treat the second term, \eqref{eq:sec2}, with an integration by parts:
\[\begin{split}
\langle \Lambda^{\sigma} E , \frac{h^{\omega,n}u^{\omega,n}}
{\rho_0}\Lambda^{\sigma} E_x \rangle  =& \frac{1}{2\rho_0}\iint_{\T^2} \partial_x 
\big( h^{\omega,n} u^{\omega,n} (\Lambda^{\sigma}E)^2\big) \,dx\,dy \\
& - \frac{1}{2\rho_0}\iint_{\T^2} \partial_x \left( h^{\omega,n} u^{\omega,n}\right) 
(\Lambda^{\sigma}E)^2 \,dx\,dy \\
= &- \frac{1}{2\rho_0}\iint_{\T^2}  \left( h^{\omega,n}_x u^{\omega,n}
+ h^{\omega,n} u^{\omega,n}_x \right) (\Lambda^{\sigma}
E)^2 dx\,dy \end{split}\]
and now Cauchy-Schwarz and the Sobolev imbedding theorem yield
\[
\left| \langle \Lambda^{\sigma} E , \frac{h^{\omega,n}u^{\omega,n}}
{\rho_0}\Lambda^{\sigma} E_x \rangle\right|  \leq C n^{-1} \|E\|_{\sigma}^2\,,
\]
where $C$ depends on $\rho_0$.
Treating the remaining terms in \eqref{eq:gr2} in the same way gives
\begin{equation} \left|
\langle \Lambda^{\sigma} \left( \diag(A_1(U^{\omega,n}))\mathcal{E}_x + 
\diag(B_1(U^{\omega,n}))\mathcal{E}_y \right), \Lambda^{\sigma}
 \mathcal{E} \rangle \right|
 \leq C n^{\max\{-1,\sigma-s+1\}} \| \mathcal{E}\|_{\sigma}^{2}\,, \label{gr2est}
\end{equation}
where the constant $C$ depends only on $\rho_0$, $\gamma$ and $h_0$.

We group the terms in \eqref{eq:gr3} to take advantage of the symmetry. 
With a change of sign we have 
\begin{multline*}
 \langle \Lambda^{\sigma} E , \Lambda^{\sigma} (h^{\omega,n}F_x)\rangle +  \langle 
 \Lambda^{\sigma} F , \Lambda^{\sigma} (h^{\omega,n}E_x)\rangle  
 \\
+ \langle \Lambda^{\sigma} E , \Lambda^{\sigma} (h^{\omega,n}G_y)\rangle +  \langle 
\Lambda^{\sigma} G , \Lambda^{\sigma} (h^{\omega,n}E_y)\rangle \\
+\langle \Lambda^{\sigma} F , \Lambda^{\sigma} (\rho_0H_x)\rangle +  \langle 
\Lambda^{\sigma} H , \Lambda^{\sigma} (\rho_0F_x)\rangle \\
+ \langle \Lambda^{\sigma} G , \Lambda^{\sigma} (\rho_0 H_y)\rangle +  \langle 
\Lambda^{\sigma} H , \Lambda^{\sigma} (\rho_0 G_y)\rangle 
\,.
\end{multline*}
Since  all the pairs are handled in the same way, we show only how to  
bound the first pair, which we rewrite using commutators as
\begin{align}
 \langle \Lambda^{\sigma} E , \Lambda^{\sigma} (h^{\omega,n}F_x)\rangle +  \langle 
 \Lambda^{\sigma} F , \Lambda^{\sigma} (h^{\omega,n}E_x)\rangle & =  \langle 
 \Lambda^{\sigma} E , \left[ \Lambda^{\sigma},  h^{\omega,n}\right] F_x\rangle + \langle 
 \Lambda^{\sigma} E ,  h^{\omega,n}\Lambda^{\sigma}F_x\rangle \label{eq:double1} \\
& + \langle \Lambda^{\sigma} F , \left[ \Lambda^{\sigma} , h^{\omega,n}\right] E_x\rangle 
+ \langle \Lambda^{\sigma} F ,  h^{\omega,n}\Lambda^{\sigma} E_x\rangle 
\label{eq:double2}
\end{align}
The first terms on the right hand side in both \eqref{eq:double1} and
 \eqref{eq:double2} are bounded by $\|\tilde h^{\omega,n}\|_{s} \|E\|_{\sigma}  
 \|F\|_{\sigma}$ from the commutator estimate \eqref{eq:est_com}. 
 We combine the second terms in \eqref{eq:double1} and \eqref{eq:double2}:
\begin{align}
 \langle \Lambda^{\sigma} E ,  h^{\omega,n}\Lambda^{\sigma}F_x\rangle + \langle 
 \Lambda^{\sigma} F ,  h^{\omega,n}\Lambda^{\sigma} E_x\rangle 
 &= \iint_{\T^2} h^{\omega,n} 
 \partial_x (\Lambda^{\sigma}E \Lambda^{\sigma} F) dx\,dy \\
& = \iint_{\T^2} \partial_x (h^{\omega,n}\Lambda^{\sigma}E \Lambda^{\sigma} F) dx\,dy -
\iint_{\T^2}  h^{\omega,n}_x\Lambda^{\sigma}E \Lambda^{\sigma} F dx\,dy\,.
\label{eq:hxEF}
\end{align}
The first term in \eqref{eq:hxEF} vanishes and the second term is estimated by $\|
\partial_x h^{\omega,n}\|_{\infty} \|E\|_{\sigma}  \|F\|_{\sigma}$ using Cauchy-Schwarz. 
Since $\|\tilde h^{\omega,n}\|_s=n^{-s}<n^{1-2s}=\|\partial_x h^{\omega,n}\|_\infty$, 
then for \eqref{eq:gr3} we have
\begin{align}
\bigg|\langle \Lambda^{\sigma} \mathcal{E} , \Lambda^{\sigma} \big( A_R(U^{\omega,n})
\mathcal{E}_x + B_R(U^{\omega,n})\mathcal{E}_y \big)\rangle
\bigg| \leq C n^{1-2s} \| \mathcal{E}\|_{\sigma}^2\,, \label{gr3est}
\end{align}
where $C$ depends only on $\rho_0$.
Note that for $n\gg 1$, $s>2$  and $1<\sigma<s-1$ we have $n^{\sigma-s+1}>n^{1-2s}$ 
and so this contribution is dominated by by the estimates \eqref{gr1est} and \eqref{gr2est}
and can be ignored.

For \eqref{eq:gr4} we use Cauchy-Schwarz and Lemma \ref{res_est} to get
\begin{equation} \left|
\left\langle\Lambda^{\sigma} H , \Lambda^{\sigma} \left( \frac{\rho_0}
{(\gamma-1)h^{\omega,n}} R_4 \right)
\right\rangle\right|  \leq C \|R_4\|_{\sigma} \| H \|_{\sigma} 
 \leq  C n^{2\sigma-3s+1} \| H \|_{\sigma}\,, \label{gr4est}
\end{equation}
where $C$ depends only on $\rho_0$, $\gamma$ and $h_0$.

Combining the estimates  \eqref{gr1est}, \eqref{gr2est} and
 \eqref{gr4est} for \eqref{eq:gr1} - \eqref{eq:gr4}, we obtain \eqref{eq:A0Et}. 

Next we use a standard treatment of symmetrizable hyperbolic systems: 
We replace the $L^2$ inner product by $\langle w, A_0(U^{\omega,n}) w\rangle$; 
 this defines an equivalent $L^2$-norm  since $A_0(U^{\omega,n})$ 
is symmetric and, for large $n$, $A_0(U^{\omega,n})\geq \kappa I>0$ with
$$\kappa=\min \left\{\rho_0, 
\frac{h_0}{2\rho_0} , \frac{\rho_0}{2(\gamma-1)h_0}\right\}\,.$$
We have
\begin{align}
 \frac{d}{dt} \| \mathcal{E} \|_{\sigma}^2  &=\frac{d}{dt}\langle \Lambda^{\sigma} 
 \mathcal{E}, A_0(U^{\omega,n}) \Lambda^{\sigma} \mathcal{E} \rangle\notag \\
& = 2\langle \Lambda^{\sigma} \mathcal{E}_t, A_0(U^{\omega,n}) \Lambda^{\sigma} 
\mathcal{E} \rangle \label{eq:AoL1}\\
&  \ \ + \langle \Lambda^{\sigma} \mathcal{E}, \big(A_0(U^{\omega,n})\big)' \Lambda^{\sigma} 
\mathcal{E} \rangle \label{eq:AoL2}
\end{align}
We write \eqref{eq:AoL1} using the symmetry of $A_0$ and a commutator as
\begin{align}
2\langle \Lambda^{\sigma} \mathcal{E}_t, A_0(U^{\omega,n}) \Lambda^{\sigma} 
\mathcal{E} \rangle & = - 2\langle \Lambda^{\sigma} \mathcal{E} ,
 [\Lambda^{\sigma} , A_0(U^{\omega,n})] 
\mathcal{E}_t \rangle \label{eq:AoL3}\\
& \ + 2\langle \Lambda^{\sigma} \mathcal{E} ,
\Lambda^{\sigma} \big( A_0(U^{\omega,n}) \mathcal{E}_t\big)  
\rangle \label{eq:AoL4}
\end{align}
The term \eqref{eq:AoL4} is estimated in \eqref{eq:A0Et}. 
For \eqref{eq:AoL3}, since $\rho_0$ is a constant and $h^{\omega,n}=h_0
+\tilde h^{\omega,n}$, we have 
\begin{align}
\langle\Lambda^{\sigma} 
\mathcal{E} , [\Lambda^{\sigma} , A_0(U^{\omega,n})] \mathcal{E}_t  \rangle 
&= \frac{1}{\rho_0}\langle\Lambda^{\sigma} E , [\Lambda^{\sigma}, \tilde h^{\omega,n}] E_t
\rangle \label{eq:comE} \\
& \ + \frac{\rho_0}{\gamma-1}\langle\Lambda^{\sigma} H ,
 [\Lambda^{\sigma}, \frac{1}{h^{\omega,n}}] H_t
\rangle \label{eq:comH}
\end{align}
By Cauchy-Schwarz and
the commutator estimate \eqref{eq:est_com}, the right hand side of \eqref{eq:comE}
 is bounded by $\|\tilde h^{\omega,n}\|_{\sigma} \|E_t\|_{\sigma-1}\|E\|_{\sigma}$ up to 
 a constant depending on $\rho_0$, and 
 in the same way \eqref{eq:comH} is bounded by 
 $ \left\| \frac{1}{h^{\omega,n}}\right\|_{\sigma} \|H_t\|_{\sigma-1}\|H\|_{\sigma}$ 
 up to a constant depending on $\gamma$ and $\rho_0$. 
 From the equation for the error \eqref{eq:errsys} we have
\begin{equation} \label{errEH}
\begin{split} 
E_t&=-\left( u^{\omega,n}E_x + \rho_0F_x + v^{\omega,n} E_y + \rho_0 G_y 
+ (u_x+v_y)E+\rho_x F + \rho_y G\right)\\
 H_t&=-(\gamma-1)\left( h^{\omega,n} (F_x+ G_y)  +(u_x+v_y)H\right) - R_4 
- u^{\omega,n}H_x - v^{\omega,n}H_y -\tilde h_x F - \tilde h_y G \,.
\end{split}\end{equation}
We cannot use the algebra property of Sobolev spaces here since 
$\sigma-1$ is not necessarily greater than $1$. 
Instead we use the following argument, which we detail here for 
$\| u^{\omega,n}E_x\|_{\sigma-1}$, 
on each term.
\begin{align*}
\| u^{\omega,n}E_x\|_{\sigma-1}^2 & = 
\|\Lambda^{\sigma-1}(u^{\omega,n}E_x)\|_{L^2}^2 
 = \| [ \Lambda^{\sigma-1}, u^{\omega,n}] E_x - u^{\omega,n} \Lambda^{\sigma-1} 
 E_x\|_{L^2}^2 \\
& \leq \| [ \Lambda^{\sigma-1}, u^{\omega,n}] E_x\|_{L^2}^2 
+\| u^{\omega,n} \Lambda^{\sigma-1}  E_x\|_{L^2}^2\,.
\end{align*}
Using the commutator estimate \eqref{eq:est_com} and the 
Sobolev embedding theorem we have
\begin{equation}
\| u^{\omega,n}E_x\|_{\sigma-1} \leq C 
( \|u^{\omega,n}\|_{\sigma+1} \|E\|_{\sigma}+  \|u^{\omega,n}\|_{s-1} \|E\|_{\sigma} ).
\end{equation}
Then the solution size estimate \eqref{eq:solsize} and the bound
\eqref{eq:appest} on the approximate solutions give
\begin{equation}
\| u^{\omega,n}E_x\|_{\sigma-1} \leq C (n^{\sigma-s+1}
+ n^{-1} )\|E\|_{\sigma} \leq C n^{\max\{-1,\sigma-s+1\}} \|E\|_{\sigma}.
\end{equation}

All the terms that arise in computing $\| E_t \|_{\sigma-1}$ and $\| H_t \|_{\sigma-1}$ 
from the right hand side of \eqref{errEH} are estimated in a similar way.
In dealing with \eqref{eq:comH}, in order to get an estimate that contains the
correct order of decay with $n$ we must replace the expressions involving 
$h^{\omega,n}$ in \eqref{errEH} with expressions in $\tilde h^{\omega,n}$,
and this can be done since we have
\[  [\Lambda^{\sigma}, \frac{1}{h^{\omega,n}}]  \big( h^{\omega,n} (F_x+ G_y)  \big) = 
\Lambda^{\sigma} (F_x+G_y)-\frac{1}{h^{\omega,n}}\Lambda^{\sigma}
 \big( h^{\omega,n} (F_x+G_y)\big)
\]
\[ =  \Lambda^{\sigma}(F_x+G_y) - \frac{h_0}{h^{\omega,n}} \Lambda^{\sigma}(F_x+G_y) - 
\frac{1}{h^{\omega,n}} \Lambda^{\sigma} \big(\tilde{h}^{\omega,n}(F_x+G_y)\big)
\]
\[ = \frac{ \tilde{h}^{\omega,n}}{h^{\omega,n}}\Lambda^{\sigma}(F_x+G_y) - \frac{1}
{h^{\omega,n}}\Lambda^\sigma \big(\tilde{h}^{\omega,n}(F_x+G_y)\big)\,.
\]
Thus,  from 
\eqref{eq:comE} and \eqref{eq:comH}  we obtain the following estimate for 
the right hand side of \eqref{eq:AoL3}:
\begin{align}\big|
\langle \Lambda^{\sigma} \mathcal{E} ,
 [\Lambda^{\sigma} , A_0(U^{\omega,n})] \mathcal{E}_t
\rangle \big|\leq C n^{\max\{-1,\sigma-s+1\} }
 \| \mathcal{E}\|_{\sigma}^2 \,,\label{eq:AoL11}
\end{align}
where $C$ depends only on $\rho_0$, $h_0$ and $\gamma$.

For \eqref{eq:AoL2} we have
\begin{align*}
\frac{d}{dt} A_0(U^{\omega,n})= \left( \begin{array}{cccc}
\frac{h^{\omega,n}_t}{\rho_0} & 0 &0 &0 \\
0& 0 &0&0\\
0& 0&0&0\\
0& 0&0&-\frac{\rho_0 h^{\omega,n}_t }{(\gamma-1)(h^{\omega,n})^2}\\
\end{array} \right),
\end{align*}
hence
\begin{align}
\langle \Lambda^{\sigma} \mathcal{E}, \big(A_0(U^{\omega,n})\big)' \Lambda^{\sigma} 
\mathcal{E} \rangle = \langle \Lambda^{\sigma} E, \frac{1}{\rho_0} \tilde h^{\omega,n}_t 
\Lambda^{\sigma} E \rangle - \langle \Lambda^{\sigma} H, \frac{\rho_0 \tilde h^{\omega,n}_t}
{(\gamma-1)(h^{\omega,n})^2} \Lambda^{\sigma} H \rangle \label{eq:A0t}\,.
\end{align}
By the definition of approximate solutions \eqref{eq:happ} we have
$\|h^{\omega,n}_t \|_{s}\leq C n^{-s}$, 
where $C$ is a constant. 
Using this last estimate in  \eqref{eq:A0t} gives
\begin{align}\big|
\langle \Lambda^{\sigma} \mathcal{E}, \big(A_0(U^{\omega,n})\big)' \Lambda^{\sigma} 
\mathcal{E} \rangle\big| \leq C n^{-s} \| \mathcal{E} \|_{\sigma}^2 \label{eq:AoL22}
\end{align}
for \eqref{eq:AoL2},
where the constant $C$ depends only on $\gamma$, $h_0$ and $\rho_0$.

Since $-s<-s+\sigma+1$, the quantity in \eqref{eq:AoL2} is dominated by \eqref{eq:AoL1},
which we have estimated in \eqref{eq:A0Et} and \eqref{eq:AoL11}.
Combining  \eqref{eq:A0Et} and \eqref{eq:AoL22} with \eqref{eq:AoL11} 
we get the bound
\begin{align}
 \frac{d}{dt} \| \mathcal{E} \|_{\sigma}^2 & \leq C \left( n^{\max\{-1,\sigma-s+1\}} \| 
 \mathcal{E}\|_{\sigma}^2 
 +  n^{2\sigma-3s+1} \| \mathcal{E} \|_{\sigma}\right)
 \end{align}
where $C$ depends only on $\rho_0$, $h_0$ and $\gamma$. 
The estimate \eqref{eq:errm} now follows by Gronwall's inequality.
\proofend

\section{Proof of Theorem \ref{th:nonuni}} \label{thmpf}

Let us now consider the two sequences of solutions $U_{1,n}(x,y,t)$ and 
$U_{-1,n}(x,y,t)$ for the initial data $U^{1,n}(x,y,0)$ and $U^{-1,n}(x,y,0)$ respectively. 
At time $t=0$ we have
\begin{align}
\| U_{1,n}(0)-U_{-1,n}(0)\|_{s}=Cn^{-1} \rightarrow 0 \mbox{ as } n \rightarrow \infty.
\end{align}
For $t>0$, by the triangle inequality we have 
\begin{equation}\begin{split}
\| U_{1,n}(t)-U_{-1,n}(t)\|_{s} 
& \geq \| U^{1,n}(t)-U^{-1,n}(t)\|_{s} - \| U^{1,n}(t)-U_{1,n}(t)\|_{s} 
 - \| U^{-1,n}(t)-U_{-1,n}(t)\|_{s} \\ \label{eq:lower1}
& \geq \| U^{1,n}(t)-U^{-1,n}(t)\|_{s} - C \|\mathcal{E}\|_{s}. 
\end{split}\end{equation}
To complete the proof, which proceeds by showing that $\|\mathcal E\|_s
\to 0$ and so we can bound the difference in actual solutions by the
difference in the approximate solutions, we need the following result.
\begin{lemma} \label{Ebound}
For $\tau\in (s,\lfloor s\rfloor +1]$ and a constant $C$ that depends on $\tau$
but not on $n$, we have
\[
\|U_{\omega,n}(t)\|_{\tau} 
  \leq C n^{\tau -s}\,,
\]
for all $t\in[0,T]$.
\end{lemma}
\proofbegin
The solution size estimate \eqref{eq:solsize} gives $\|U(t)\|_\tau\leq 
C(\tau, d)\|U(0)\|_\tau$ for all data with $\|U(0)\|_\tau\leq d$, for any $\tau>2$
for which $\|U(0)\|_\tau$ is defined, and for all $t\in[0,T)$ where $T$ also
depends on $d$ and on $\tau$.
 Furthermore (see Corollary 2 to Theorem 2.2 in Majda \cite{Maj2}), if $T$ is
 a maximum lifespan, then either $U$ leaves every compact subset of $G$
 (the subset of phase space in which the system is symmetrizable hyperbolic)
 or $\|\nabla U(t)\|_{L^\infty} + \|U_t(t)\|_{L^\infty}\to\infty$ as $t\to T$.
 This means, for our solutions, since the data are in $H^\tau$ for all $\tau >0$
 and we assume we have identified a $T<T_{\textrm{crit} }$, where 
 $T_{\textrm{crit} }$ is the value beyond which a solution in $H^s$ no
 longer exists,  that the solution remains in $H^\tau$ for $t<T$ and any 
 $\tau>s$.
 (Here we note that $s>2$ so $U$ and its first derivatives are bounded,
 both pointwise and in $H^s$, for $t\in[0,T]$.)
 
 However, in the estimate on the solution size \eqref{eq:solsize}, 
 the constant $C$ depends on
 $d=\|U(0)\|_\tau$, and this is bounded (by unity) only for $\tau\leq s$.
 If $\tau>s$, then $\|U(0)\|_\tau \to \infty$ with $n$.
 To use the interpolation result, \eqref{interpol} below, we need to apply \eqref{eq:solsize},
 with a constant independent of $n$, for some
 value of $\tau>s$.
  
 We obtain a bound for  $\tau= \lfloor s \rfloor +1$, where $\lfloor s \rfloor$ is
 the greatest integer in $s$, as follows. 
 Let $\alpha$ with $|\alpha|=\tau$ be a
 multi-index corresponding to any $\tau^{\textrm{th}}$ order derivative.
 There are $\tau +1$ such derivatives; define $V_i= D^{(\tau+1-i, i-1)}U$.
 Differentiating \eqref{classical} $|\alpha|$ times for all $\alpha$ with $|\alpha|=\tau$ leads to
 \[ (\partial_t + A(U) \partial_x + B(U)\partial_y)V_i+ 
 \sum_{j=1}^{\tau+1}M_j V_j+f_i(D^{\beta}U;|\beta|\leq\tau -1)=0 
 \quad \mbox{ for } \quad i=1, ... , \tau+1\,,
 \]
 where the $M_j$ are block diagonal matrices that depend only on $U$ and $DU=(U_x,U_y)$.
 Thus, $V=(V_1, V_2, ... , V_{\tau+1})$ is the solution of a linear symmetrizable 
 hyperbolic system with bounded coefficients.
 The secular term $f=(f_1, f_2, ..., f_{\tau+1})$ is also bounded, so the usual energy estimates,
 applied to the symmetrized system, yield a bound for $V$ that depends
 on the value of $V(0)$ (and as usual on $\rho_0$, $h_0$ and $\gamma$,
 and our original choice for $T$, but on nothing else).
 From \eqref{eq:initial} and \eqref{eq:appest}, a bound for $V(0)$ is $Cn^{\tau-s}$.
  This
 gives the bound stated in the Lemma for the actual solution $U_{\omega,n}$,
 for any $\tau \leq \lfloor s\rfloor +1$.
 \proofend
 
Theorem \ref{prop:error} gives a bound for $\|\mathcal{E}\|_{\sigma}$, for
$1<\sigma < s-1$. 
We use interpolation (Theorem 5.2 in \cite{Adams}) between $\sigma$ and 
$\tau=\lfloor s \rfloor +1$ to obtain a bound for $\|\mathcal{E}\|_{s}$:
\be \label{interpol} 
 \big\|\mathcal{E}\big\|_{ {s} }  \leq\big\| \mathcal{E}\big\|_ \sigma^{\alpha}  
\big\| \mathcal{E}\big\|_\tau^{\beta}, \quad \textrm{where}
\quad \alpha= \frac{ \tau -s}{\tau - \sigma}\,,\quad \beta =
\frac{  s -\sigma }{\tau - \sigma} \,.
\ee
Now,  assume we have fixed a
compact set $G_2$ with $\rho\geq \rho_0/2$, say, and once $c$
in  Theorem \ref{prop:error} is
bounded then so is $ct$ for  $t\leq T$, so
$\|\mathcal E\|_\sigma\leq Cn^\nu$, where $\nu=\max\{2\sigma -3s+2, \sigma - 2s\}$,
 and thus the exponent of $n$
in $\|\mathcal E\|_s$ is
$$ \alpha\nu + \beta(\tau - s) =\frac{(\tau-s)\nu}{\tau-\sigma} +
\frac{(s-\sigma)(\tau-s)}{\tau-\sigma}
=\frac{\tau-s}{\tau-\sigma}\big(\max\{\sigma-2s+2, -s\}\big)
$$
and this is negative since we have assumed $\sigma < s-1$ and $s>2$.
Thus, the $H^s$ error in the approximate solutions tends to zero as
$n\to\infty$, and we can estimate the difference between the actual
solutions by the difference in the approximate solutions.

Using trigonometric identities, we have
\begin{multline*}
 U^{1,n}(t)-U^{-1,n}(t)   =\\  \begin{pmatrix}
0 \\ 
\displaystyle \frac{2}{n}+\frac{1}{n^s}(\cos(ny-t)-\cos(ny+t)) \\ \\
\displaystyle  \frac{2}{n}+\frac{1}{n^s}(\cos(nx-t)-\cos(nx+t)) \\ \\
\displaystyle  \frac{1}{n^{2s}}[\sin(nx-t)\sin(ny-t)-\sin(nx+t)\cos(ny+t)] 
 \end{pmatrix}  
  =   \begin{pmatrix}
0 \\
\displaystyle \frac{2}{n}+\frac{2}{n^s}\sin ny\sin t \\ \\
\displaystyle \frac{2}{n}+\frac{2}{n^s}\sin nx\sin t \\ \\
\displaystyle -\frac{1}{n^{2s}}\sin(nx+ny)\sin 2t 
 \end{pmatrix}\,.
\end{multline*}

Then the estimate \eqref{eq:lower1} implies
\[
\liminf_{n\rightarrow \infty} \| U_{1,n}(t)-U_{-1,n}(t)\|_{s}  \geq \liminf_{n\rightarrow \infty}  
\| U^{1,n}(t)-U^{-1,n}(t)\|_{s}  \geq C \sin t\,.
\]
This completes the proof of nonuniform dependence.
 
 \section{Conclusions}
 \label{conc}
 
\begin{figure}
\centering
\includegraphics[width=.48\textwidth,clip=,viewport=130 160 485 590]{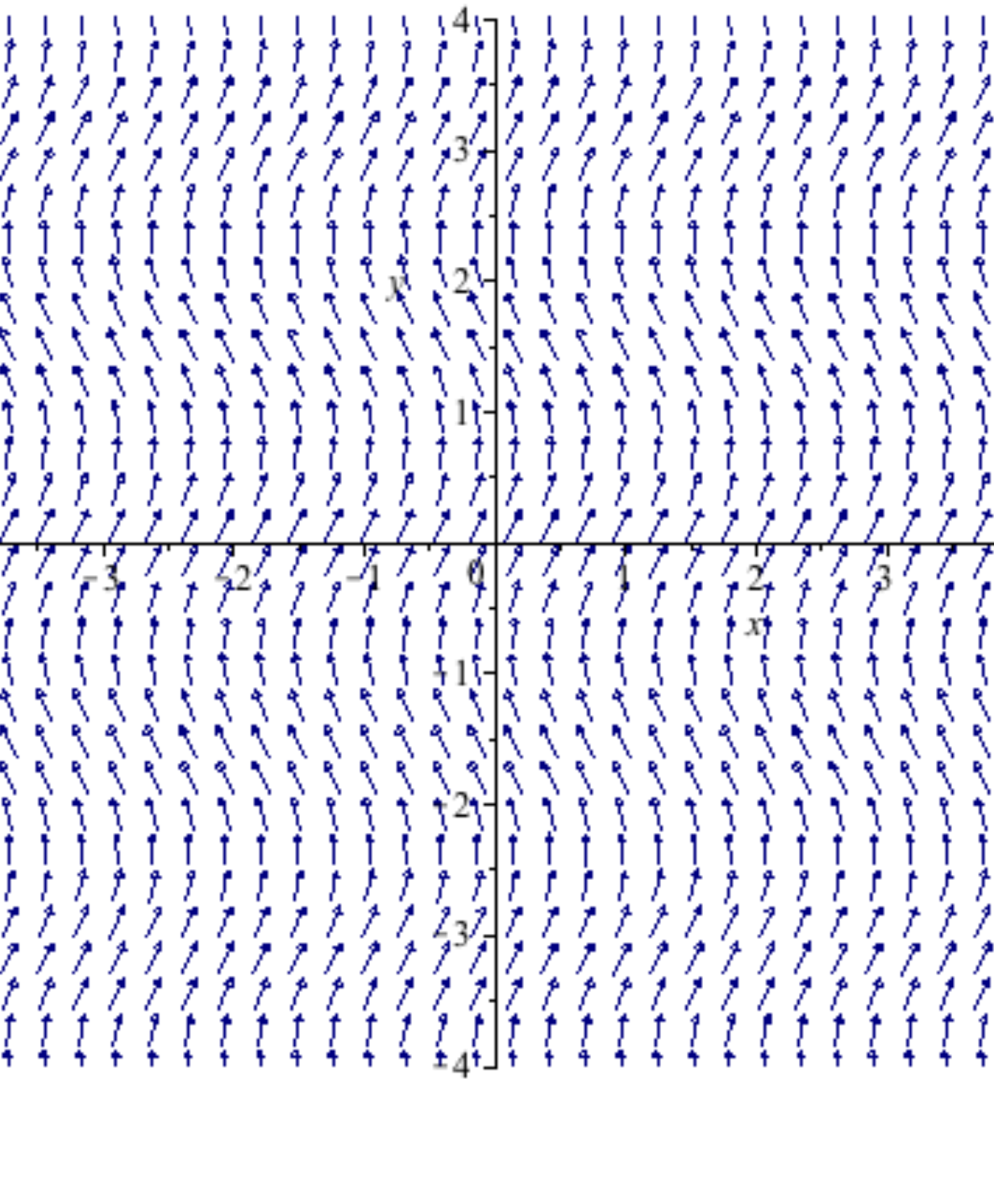} \quad \includegraphics[width=.48\textwidth,clip=,viewport=120 150 490 600]{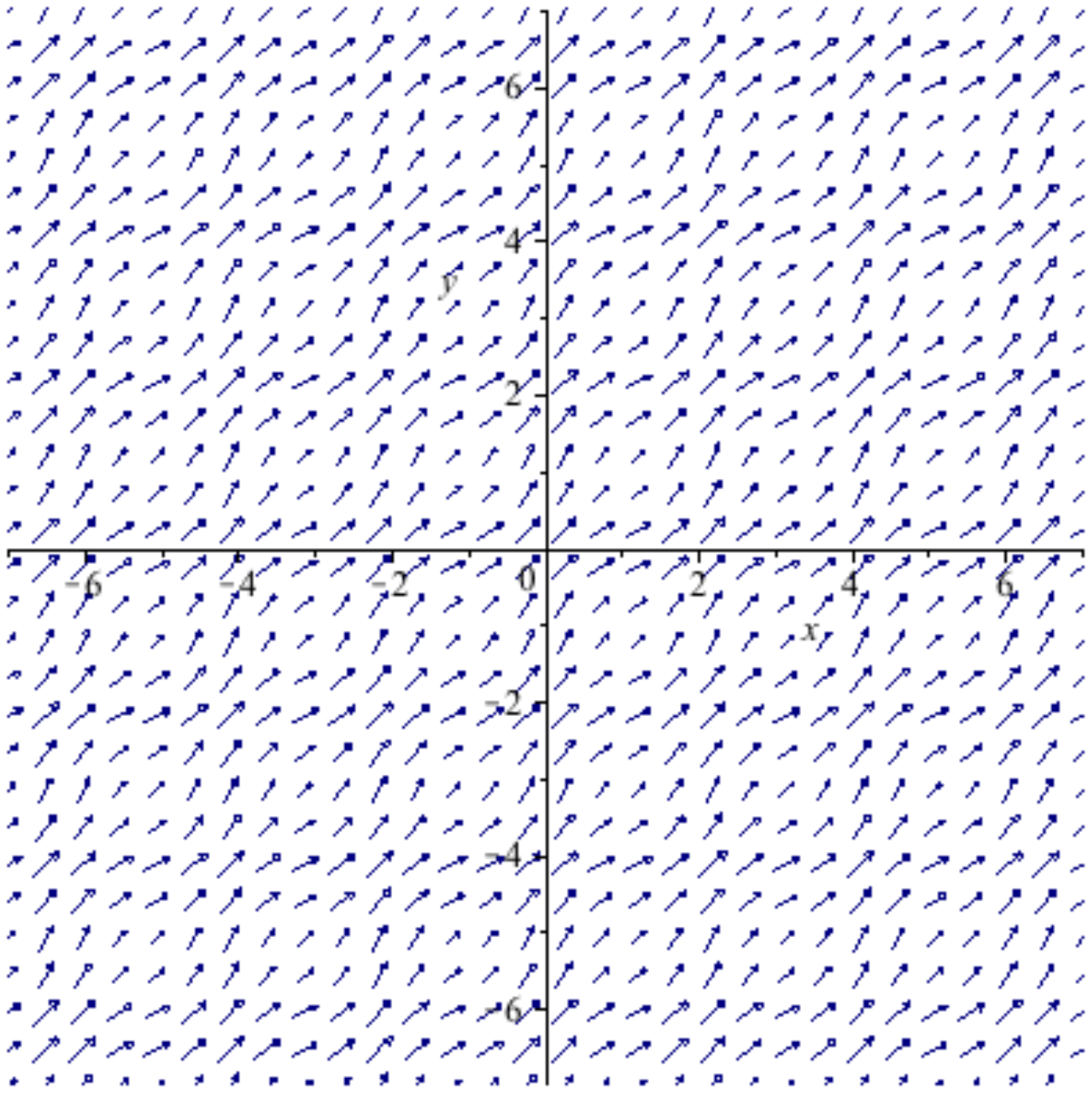}
\caption{Velocity Fields for the Constant Density (Left) and Approximate (Right)
Solutions}\label{sketch}
\end{figure}

This paper shows that periodic solutions of the compressible gas dynamics
equations in two space dimensions exhibit nonuniform dependence on
initial conditions, by a mechanism very similar to that governing the
incompressible system. 
Both the constant-density construction of Section \ref{cstden} and the approximate
solutions based on the Himonas-Misio\l ek model take an initial condition
consisting of a uniform motion with a smaller oscillatory motion superimposed
on it.
We sketch the initial velocity fields for typical members of each series in
Figure \ref{sketch}.
The constant-density and constant-pressure solution is not completely
trivial.
It is also a solution to the incompressible system, somewhat
simpler than the one devised by Himonas and Misio\l ek.
It persists for all time, without the formation of shocks.
There may be  other families of solutions and
approximate solutions with similar structure.
The actual solutions corresponding to our approximation \eqref{eq:happ}
do not have constant density or pressure.

 The conclusions to be drawn from this demonstration are of two types.
 First, ``nonuniform dependence on data'' in the sense of this paper
 can be contrasted to ``uniform dependence'' in the sense of Kato's 
 original well-posedness proof.
 Second, it is worth calling attention to the nature of the solutions we have
 constructed, as they are solutions of a hyperbolic system (compressible
 flow) that is closely related to a system that is not hyperbolic (incompressible
 flow).
 
 We look at these separately.
 
 \subsection{The Meaning of Non-Uniform Dependence}
The failure of uniform dependence on the data is instantaneous
 and is a property of classical solutions.
 It does not appear to tell us anything about properties of weak solutions
 (existence of which, for the multidimensional compressible Euler system,
 is an open problem).
 In his important paper \cite[Theorem III(b)]{Kato}, Kato
 proves {\em uniform} dependence of solutions on the data in the situation
 where a limiting initial condition in $H^s$ is approximated in $H^s$ by a
 sequence of initial conditions $\{U_0^n\}$. 
 That is not the case for our data.
 While the difference between corresponding terms in 
 our sequences $U^{1,n}$ and $U^{-1,n}$
 converges to zero in $H^s$, neither sequence alone converges in $H^s$.
 In verifying the error bounds claimed for the approximate
 solutions, one can see that the cancellation between the ``low frequency'' terms
 ($\pm 1/n$ in this case) and the high frequency oscillatory terms is 
 a result of nonlinearities in the system. 
 This creates the possibility of the nonuniformity demonstrated here.
 A similar type of cancellation,  differing in detail, is used in our companion
 paper \cite{HoKeTi} to obtain a nonuniformity result for solutions defined on
 on the plane, rather than on a torus.

 \subsection{Linear and Nonlinear Behavior in Gas Dynamics}
 
 It is also interesting to compare the nonuniform sequences of solutions 
 we have constructed here with the sequences Himonas and
 Misio\l ek \cite{Himonas-Misiolek} used in their proof of nonuniform
 dependence for the incompressible system.
 That system takes the form of three equations, for velocity and pressure:
 \begin{equation} \label{incomp} \nabla\cdot{\bf u}=0\,, \quad
{\bf u}_t+\nabla\left( {\bf u} \times{\bf u} \right)+\nabla p=0\,.\end{equation}
 This system is not hyperbolic; to the extent that its characteristics can
 be compared to those of \eqref{eq:ncompgas1}, one could say that
 the acoustic characteristics in  \eqref{eq:ncompgas1} (those associated
 with the ``speed of sound'', and also the pair that are genuinely nonlinear
 in the sense of conservation laws) have become infinite in \eqref{incomp}.
 (This is more correctly stated in terms of the Mach number -- the ratio of
the fluid velocity to  the characteristic speed.
 The system \eqref{incomp} represents a flow in which the Mach number
 has become zero.)
 
 Our exhibition of nonuniform behavior in a hyperbolic system related to
 the incompressible system indicates that the nonuniform dependence is
 \begin{itemize}
 \item[(a)] hyperbolic in nature, and
 \item[(b)] based in the linear characteristics of the hyperbolic system,
 which are shared with the incompressible system -- that is, the shear
 or entropy waves.
 \end{itemize}

Finally, we observe that a simple adaptation of
the constant-density example of Section \ref{cstden} also proves nonuniform
dependence on data 
 for the isentropic gas dynamics system -- the
system formed from the first three equations of \eqref{eq:ncompgas1} by
assuming that the pressure is a given function of the density.
That system, of course, has only a single linear family, corresponding
to shear waves.

 \bibliography{biblio}
 
 \end{document}